\documentclass[a4paper,12pt,leqno]{article}
\usepackage{latexsym}
\usepackage[all]{xy}

\usepackage{amssymb} 
\usepackage{amsmath} 
\usepackage{theorem}

\def\P{{\mathbf{P}}}
\def\Z{{\mathbb{Z}}}
\def\K{{\mathbb{K}}}
\def\CC{{\mathbb{C}}}
\def\R{{\mathbb{R}}}
\def\A{{\mathcal{A}}}
\def\B{{\mathcal{B}}}

\DeclareMathOperator{\rank}{rank}
\DeclareMathOperator{\codim}{codim}

\DeclareMathOperator{\Der}{Der}

\DeclareMathOperator{\Shi}{Shi}
\DeclareMathOperator{\Cat}{Cat}

\numberwithin{equation}{section}

\newcommand{\owari}{\hfill$\square$}

\theoremstyle{break}
\newtheorem{theorem}{Theorem}[section]
\newtheorem{prop}[theorem]{Proposition}

\newtheorem{lemma}[theorem]{Lemma}
\newtheorem{define}[theorem]{Definition}

\newtheorem{rem}[theorem]{Remark}

\newtheorem{conj}[theorem]{Conjecture}

\title{Addition-deletion theorem for free hyperplane arrangements and 
combinatorics}
\author{Takuro Abe
\footnote
{
Institute of Mathematics for Industry,
Kyushu University,
Fukuoka 819-0395, Japan.
Email:abe@imi.kyushu-u.ac.jp. 
}
}
\date{\today}

\begin{document}

\maketitle

\begin{abstract}
In the theory of hyperplane arrangements, the most important and difficult problem is 
the combinatorial dependency of several properties. In this atricle, 
we prove that Terao's celebrated addition-deletion theorem for free arrangements 
is combinatorial, i.e., 
whether you can apply it depends only on the intersection lattice of 
arrangements. The proof is based on a classical technique. 
Since some parts are already completed recently, we prove the rest 
part, i.e., the combinatoriality of the addition theorem. 
As a corollary, we can define a new class of free arrangements called the additionally free arrangement of hyperplanes, which can be constructed 
from the empty arrangement by using only the addition theorem. Then we can show that 
Terao's conjecture is true in this class. As an application, we can show that 
every ideal-Shi arrangement is additionally free, implying that 
their freeness is combinatorial.
\end{abstract}

\section{Introduction}

Let $V=\K^\ell$, $S=\mbox{Sym}(V^*) =\K[x_1,\ldots,x_\ell]=\oplus_{d=0}^\infty S_d$ with the standard polynomial degree grading, and $\Der S=\oplus_{i=1}^\ell S \partial_{x_i}$. An \textbf{arrangement of hyperplanes} $\A$ is a finite set of linear hyperplanes in 
$V$. 
In the study of hyperplane arrangements, the most interesting and difficult problem is 
to determine whether some property of $\A$ depends only on combinatorics. For example, when $\K=\CC$, the cohomology ring 
$H^*(V \setminus \cup_{H \in \A} H;\Z)$ is shown to be combinatorial in \cite{OS}. However in general, to determine whether some property of $\A$ is combinatorial or not 
is a hard problem. In this atricle, we will show that an important property of 
algebra of $\A$ is combinatorial, which completes the recent sequental works by the author 
in \cite{A}, \cite{A2}, \cite{A4} and \cite{A5}. For that, let us introduce a notation for algebra of $\A$. 

For each $H \in \A$, let $\alpha_H$ denote a fixed defining linear form of $H$. Now the \textbf{logarithmic derivation module} $D(\A)$ of $\A$ is defined by 
$$
D(\A):=\{\theta \in \Der S \mid \theta(\alpha_H) \in S \alpha_H\ (\forall H \in \A)\}.
$$
$D(\A)$ is an $S$-graded reflexive module of rank $\ell$, but not a free $S$-module in general. 
We say that $\A$ is \textbf{free} with \textbf{exponents} $\exp(\A)=(d_1,\ldots,d_\ell)$ if there is a homogeneous $S$-basis $\theta_1,\ldots,\theta_\ell$ for $D(\A)$ such that $\deg \theta_i=d_i$ for 
all $i$. Here we say that a non-zero derivation $\theta \in \Der S$ is \textbf{homogeneous of degree $d$} if 
$\theta(\alpha) \in S_d$ for all $\alpha \in V^*$. For example, the \textbf{Euler derivation} $\theta_E:=\sum_{i=1}^\ell x_i \partial_{x_i}$, which is always contained in $D(\A)$, is homogeneous of degree $1$. 

Because of several important results which relate the algebra $D(\A)$, topology and 
combinatorics of $\A$, the logarithmic derivation module and free arrangements have been being studied. For example, see Theorem \ref{factorization}. However, it is not easy to check the freeness. The most useful result to check the freeness is the following addition-deletion theorem due to Terao in 1980.

\begin{theorem}[Addition-deletion theorem, \cite{T1}]
Let $H \in \A$, $\A':=\A \setminus\{H\}$ and $\A^H:=\{L \cap H \mid L \in \A'\}$. Then two of the following three imply the third:

\begin{itemize}
\item[(1)]
$\A$ is free with $\exp(\A)=(d_1,\ldots,d_{\ell-1},d_\ell)$.

\item[(2)]
$\A'$ is free with $\exp(\A')=(d_1,\ldots,d_{\ell-1},d_\ell-1)$.

\item[(3)]
$\A^H$ is free with $\exp(\A^H)=(d_1,\ldots,d_{\ell-1})$.
\end{itemize}

In particular, all the three above are true if both $\A$ and $\A'$ are free.
\label{additiondeletion}
\end{theorem}

Though Theorem \ref{additiondeletion} is proved almost 40 years before, it is the most useful to construct free arrangements, and to check the freeness and non-freeness even now. Related to this the relation between the freeness and combinatorics of $\A$ is well-studied too. Namely, the \textbf{intersection lattice} $L(\A)$ is defined by 
$$L(\A):=\{
\cap_{H\in \B} H \mid 
\B
\subset 
\A\}.
$$
This is considered to be combinatorial information of $\A$. Then \textbf{Terao's conjecture} asks whether 
the freeness depends only on $L(\A)$. Recently, the relation between Theorem \ref{additiondeletion} and combinatorics is intensively studied. See Theorems \ref{division} and  \ref{deletioncomb} 
for example. In particular, the deletion theorem ((1)+(3) $\Rightarrow$ (2) in 
Theorem \ref{additiondeletion}) is shown to be combinatorial in Theorem \ref{deletioncomb}. 
Actually, in these results, not the total structure of $L(\A)$ but the following three 
combinatorial objects play the key roles:

(i)\,\, The \textbf{localization} 
$$\A_X:=\{H \in \A \mid X \subset H\}$$ 
of $\A$ at $X \in L(\A)$.

(ii)\,\,
The \textbf{restriction} $$
\A^H:=\{L \cap H \mid L \in \A \setminus \{H\}\}
$$
of $\A$ onto $H$, 
which is an arrangement in $H =\K^{\ell-1}$. Also, for $X \in L(\A^H)$, define
$$
\A_X^H:=(\A_X)^H=(\A^H)_X.
$$
It is easy to show the second equality above.

(iii)\,\, 
The \textbf{characteristic polynomial} $\chi(\A;t)$ of $\A$, see Definition \ref{cha}.

To state the main result in this article, let us introduce the following new 
combinatorial property:

\begin{define}
Let $H \in \A$. We say that $\A$ is \textbf{divisional along $H$} if 
$\chi(\A_X^H;t) \mid \chi(\A_X;t)$ for all $X \in L(\A^H)$. 
We say that $\A$ is \textbf{locally (resp. globally) divisional along $H$} if 
$\chi(\A_X^H;t) \mid \chi(\A_X;t)$ for all $X \in L(\A^H) \setminus \{\cap_{Y \in 
\A^H} Y\}$ (resp. $\chi(\A^H;t) \mid \chi(\A;t))$.
\label{divisional}
\end{define} 

By the divisibility along $H$, we can prove a combinatorial version of Theorem \ref{additiondeletion} as follows:

\begin{theorem}[Combinatorial addition-deletion theorem]
Let $H \in \A$, $\A':=\A \setminus\{H\}$ and $\A^H:=\{L \cap H \mid L \in \A'\}$. Then two of the following four imply the third and fourth:

\begin{itemize}
\item[(i)]
$\A$ is free.

\item[(ii)]
$\A'$ is free.

\item[(iii)]
$\A^H$ is free, and $\A$ is globally divisional along $H$.

\item[(iv)]
$\A$ is divisional along $H$.
\end{itemize}

In particular, all the four above are true if (iii) is true.
\label{additiondeletion2}
\end{theorem}

In Theorem \ref{additiondeletion2}, $(iii) \Rightarrow (i), (ii), (iv)$ is 
the division theorem in \cite{A2} (Theorem \ref{division}). $(i)+(ii) \Rightarrow 
(iii)+(iv)$ follows from Theorem \ref{additiondeletion} and Lemma \ref{LD}. 
$(i)+(iv) \Rightarrow 
(ii)+(iii)$ follows from \cite{A4} (Theorem \ref{deletioncomb}). So the rest part is 
$(ii)+(iv) \Rightarrow 
(i)+(iii)$, i.e., whether the addition theorem ((2)+(3) $\Rightarrow$ (1) in 
Theorem \ref{additiondeletion}) 
is combinatorial or not. 
This is the main part of 
this article, as follows.

\begin{theorem}
For $H \in \A$, let $\A':=\A \setminus \{H\}$. 
Assume that $\A'$ is free. 
Then $\A$ is free 
if and only if $\A$ is divisional along $H$. Hence whether $\A$ is free or not 
depends only on $L(\A)$ when $\A'$ is free. 
\label{main}
\end{theorem}

From the viewpoint of Theorems \ref{additiondeletion2} and \ref{main}, the addition-deletion theorem may be regarded as the result that implies the freeness of all the 
triple $\A,\A \setminus \{H\},\A^H$ by (A) the freeness of one of them, and (B) one 
combinatorial condition. So the following is the new interpretation of Theorem \ref{additiondeletion} from this viewpoint. 

\begin{theorem}[Free triple theorem]
For $H \in \A$, assume that one of the following holds:
\begin{itemize}
\item[(i)]
$\A$ is free and $\A$ is divisional along $H$.
\item[(ii)]
$\A \setminus \{H\}$ is free and $\A$ is divisional along $H$.
\item[(iii)]
$\A^H$ is free and $\A$ is 
globally divisional along $H$.
\end{itemize}
Then all of $\A,\A \setminus \{H\},\A^H$ are free.
\label{three}
\end{theorem}

Also, we can define the following class of free arrangements, which is very natural to consider, and 
by which we have constructed a lot of free arrangements. 

\begin{define}[Additionally free arrangements]
The set $\mathcal{AF}$ consists of 
arrangements $\A$ such that there is a filtration $$
\emptyset=\A_0 \subset \A_1 \subset \cdots \subset \A_n=\A
$$
such that 
$\A_i$ is divisional along $H_i$, where $\A_i \setminus \A_{i-1}=\{H_i\}$ for 
$=1,\ldots,n$. Such a filtration is called the \textbf{additional filtration} of $\A$, and an arrangement $\A \in \mathcal{AF}$ is called an  
\textbf{additionally free arrangement}.
\label{AF}
\end{define}

\begin{rem}
The additional filtration is the same as the \textbf{free filtration} introduced in \cite{AT0}. Because we call $\mathcal{AF}$ the set of additionally free arrangements, we use the terminology ``additional filtration'' here.
\end{rem}

In other words, an arrangement is additionally free if and only if that can be constructed from the empty arrangement by using the addition theorem, which is a really natural class of free arrangements. In fact, we can show the following contribution to Terao's conjecture in terms of additionally free arrangements.

\begin{theorem}
$\A$ is free if $\A \in \mathcal{AF}$. In particular, Terao's 
conjecture is true in the class of additionally free arrangements.
\label{AFTC}
\end{theorem}

It was known that Terao's conjecture is true 
in the class of inductively free arrangements (Definition \ref{IF}) which can be 
constructed by using the addition and restriction theorem. In \cite{A2}, the new class so called the 
divisionally free arrangement is defined (Definition \ref{DF}) in which Terao's conjecture is 
true, and can be constructed by only using the division theorem (Theorem \ref{division}). Additionally free arrangements are the other direction of this generalization. Actually, we can show that the freeness of ideal-Shi arrangements depends only on 
combinatorics  by using the additional freeness (Theorem \ref{idealShi}). Moreover, we define the largest class of free arrangements $\mathcal{SF}$ in which Terao's conjecture is true. 
The class $\mathcal{SF}$ is called the \textbf{stair-free arrangements}, constructed in \S4 by combining both additionally and 
divisionally free arrangements. 

The organization of this article is as follows. In \S2 we recall several results on arrangements. 
In \S3 we prove Theorem \ref{main}. In \S4 we study additionally free arrangements and its extension joining with the divisionally free arrangements. In \S5 we apply our main result to 
ideal-Shi arrangements, proving that their freeness depends only on combinatorics. 
\medskip

\noindent
\textbf{Acknowledgements}. 
The author is 
partially supported by KAKENHI, JSPS Grant-in-Aid for Scientific Research (B) 16H03924, and 
Grant-in-Aid for Exploratory Research 16K13744. 

\section{Preliminaries}

In this section let us collect several definitions and results on arrangements which will be used in this article. See \cite{OT} for general reference. First let us recall several combinatorics of arrangements.

\begin{define}
For the intersection lattice $$
L(\A):=\{\cap_{H \in \B} H \mid \B \subset \A\}
$$
of an arrangement $\A$, 
define 
$$
L_k(\A):=\{X \in L(\A) \mid \codim_V X=k\}.
$$
The \textbf{rank $r(\A)$} of an arrangement $\A$ is defined by 
$$
r(\A):=\codim \cap_{H \in \A}H.
$$
For $X \in L(\A)$, the \textbf{localization} 
$\A_X$ of $\A$ at $X$ is defined as 
$$
\A_X:=\{H \in \A \mid X \subset H\}.$$
For $H \in \A$, the \textbf{restriction} 
$\A_X$ of $\A$ onto $H$ is defined as 
$$
\A^H:=\{H \cap L\mid L\in \A \setminus \{H\}\}.$$
Let $$
0_X:=\cap_{H \in \A_X}H \in L(\A_X)$$
for $X \in L(\A)$. So 
$\rank (\A_X)=\codim_V 0_X$. 
\end{define}

Next recall the most important combinatorial invariant of $\A$. 

\begin{define}
(1)\,\,
The \textbf{M\"obius function} $\mu:L(\A) \rightarrow \Z$ is defined by 
$\mu(V)=1$, and by 
$$
\mu(X)=-
\sum_{Y \in L(\A),\ X \subsetneq 
Y \subset V} \mu(Y)
$$
for $Y \in L(\A) \setminus \{V\}$. 

(2)\,\, 
The \textbf{characteristic polynomial} $\chi(\A;t)$ is defined by 
$$
\chi(\A;t):=\sum_{X \in L(\A)} \mu(X)t^{\dim X}.
$$
\label{cha}
\end{define}

\begin{rem}
By definition $\chi(\A;t)$ is a combinatorial invariant of $\A$. Moreover, when $\K=\mathbb{C}$, 
$(-t)^\ell \chi(\A;-1/t)$ coincides with the topological Poincar\'e polynomial of 
$\mathbb{C}^\ell \setminus \cup_{H \in \A} H$. So it is a topological invariant too.
\end{rem}

To compute $\chi(\A;t)$, the following \textbf{deletion-restriction formula} is the most useful. For the proof, 
see Corollary 2.57 in \cite{OT} for example.

\begin{theorem}
For $H \in \A$, it holds that 
$$
\chi(\A;t)=\chi(\A \setminus \{H\};t)-\chi(\A^H;t).
$$
\label{DR}
\end{theorem}

We also use the following property, which directly follows from the definition of M\"obius 
function in Definition \ref{cha} (1). 

\begin{prop}
When $\A$ is not empty, 
$\chi(\A;t)$ is divisible by $t-1$. Define the 
\textbf{reduced characteristic polynomial} 
$\chi_0(\A;t):=\chi(\A;t)/(t-1)$. 
\label{reducedcha}
\end{prop}

The localization behaves nicely in the sense of freeness as follows, see 
Theorem 4.37 in \cite{OT} for the proof.

\begin{theorem}
If $\A$ is free, then $\A_X$ is free 
for all $X \in L(\A)$.
\label{localization}
\end{theorem}

We will see the relation between the localization and freeness for details in the following.

\begin{define}
Let $X \in L_k(\A)$, and let 
$\emptyset_\ell$ denote the \textbf{empty arrangement} in $\K^\ell$. Then 
$\A_X$ decomposes into 
$$
\A_X=\B_X \times \emptyset_{\ell-k}$$
for some arrangement $\B_X$ in $V/X$ with 
$\mbox{rank}(\B_X)=k$, and 
the empty arrangement $\emptyset_{\ell-k}$ in $X$. 
We say that $\B_X$ is the \textbf{essential part} of $\A_X$.
\label{essentialpart}
\end{define}

For the proof of the following, see Proposition 4.14 in \cite{OT} for 
example.

\begin{prop}
Let $X \in L_k(\A)$ and $\B_X$ be the essential part of $\A_X$. Then 
$\A_X$ is free if and only if $\B_X$ is free.
\label{directfree}
\end{prop}

Let us introduce the local version of freeness.

\begin{define}
We say that $\A$ is \textbf{locally free} if $\A_X$ is free for all $X \in L(\A) \setminus \{
0_V\}$. 
\label{locallyfree}
\end{define}

It is known that $\A$ is locally free if and only if $\widetilde{D(\A)}$ is a vector bundle on 
$\P^{\ell-1}=\mbox{Proj}(S)$. 
We use the following famous criterion for freeness. For the proof of the arrangement version, see Theorem 4.19 in \cite{OT}. 

\begin{theorem}[Saito's criterion, \cite{Sa}]
$\A$ is free if and only if there are homogeneous derivations $\theta_1,\ldots,\theta_\ell \in D(\A)$ such that they are independent over $S$, and that 
$|\A|=\sum_{i=1}^\ell \deg \theta_i$.
\label{Saito}
\end{theorem}

The freeness is related to combinatorics and topology of arrangements by the following 
famous result.

\begin{theorem}[Teraro's factorization, Main Theorem in \cite{T2}]
Let $\A$ be free with $\exp(\A)=(d_1,\ldots,d_\ell)$. Then 
$$
\chi(\A;t)=\prod_{i=1}^\ell (t-d_i).
$$
\label{factorization}
\end{theorem}

Based on the division property in Definition \ref{divisional} with the theory of multiarrangements, following two generalizations of 
Theorem \ref{additiondeletion} have been 
obtained. 

\begin{theorem}[Division theorem, Theorem 1.1, \cite{A2}]
Let $H \in \A$. If $\A^H$ is free, and $\chi(\A^H;t) \mid 
\chi(\A;t)$, then $\A$ is free.
\label{division}
\end{theorem}

\begin{theorem}[Theorem 1.2, \cite{A4}]
Let $\A$ be free, $H \in \A$ and let 
$\A':=\A \setminus \{H\}$. Then $\A'$ is free if and only if $\A$ is divisional along $H$.
\label{deletioncomb}
\end{theorem}

There are several class of free arrangements in which Terao's conjecture is true. Let us 
recall two of them. 

\begin{define}[\cite{OT}]
Let $\mathcal{IF}_\ell$ denote the set of arrangements in $\K^\ell$ defined in the following manner. 

\begin{itemize}
\item[(1)] $\emptyset_\ell \in \mathcal{IF}_\ell$ for all $\ell \ge 0$. 

\item[(2)]
$\A \in \mathcal{IF}_\ell$ if there is $H \in \A$ such that $\A \setminus \{H\} 
\in \mathcal{IF}_\ell$, and $\A^H \in \mathcal{IF}_{\ell-1}$. 
\end{itemize}
The arrangement $\A$ belonging to 
$$
\mathcal{IF}:=\cup_{\ell \ge 0} \mathcal{IF}_\ell
$$
is called the \textbf{inductively free arrangement}. 
\label{IF}
\end{define}

\begin{define}[\cite{A2}]
The class $\mathcal{DF}_\ell$ of arrangements in $\K^\ell$ consists of arrangements 
$\A$ such that, there is $X_i \in L_i(\A)\ (i=0,\ldots,\ell-2)$ such that 
$$
X_0=V \supset X_1 \supset \cdots \supset X_{\ell-2},
$$
and $\chi(\A^{X_i};t) \mid \chi(\A^{X_{i-1}};t)$ for $i=1,\ldots,\ell-2$. The set 
$$
\mathcal{DF}:=\cup_{\ell \ge 0} \mathcal{DF}_\ell
$$
is called the set of \textbf{divisionally free arrangements}, and the above flag $\{X_i\}_{i=0}^{\ell-2}$ of 
$\A$ is called a \textbf{divisional flag} of $\A$.
\label{DF}
\end{define}

\begin{theorem}[Theorems 1.3 and 1.6, \cite{A2}]
$$
\mathcal{IF} \subsetneq \mathcal{DF},
$$
and $\A \in \mathcal{DF}$ is free. Moreover, Terao's conjecture is true for divisionally 
free arrangements.
\label{IFDF}
\end{theorem}



For the proof of our main result, let us recall the following key results. For the proof, see Proposition 4.41 in \cite{OT} for example. 

\begin{prop}[\cite{T1}]
Let $H \in \A$. Then there is a polynomial $B$ of degree $|\A|-|\A^H|-1$ such that 
$\alpha_H \nmid B$, and 
$$
\theta(\alpha_H) \in (\alpha_H,B)
$$
for all $\theta \in D(\A \setminus \{H\})$. 
\label{B}
\end{prop}

The following submodule of $D(\A)$ plays a key rote too in the proof of Theorem \ref{main}.

\begin{prop}[e.g., Lemma 1.33, \cite{Y3}]
For $H \in \A$, define 
$$
D_H(\A):=\{
\theta \in D(\A) \mid \theta(\alpha_H)=0\}.
$$
Then 
$$
D(\A)=S\theta_E \oplus D_H(\A).
$$
In other words,
$$
D(\A)/S\theta_E \simeq D_H(\A)$$
as $S$-graded modules. 
Thus $\theta_E$ is a part of basis when $\A \neq \emptyset$ is free. So $\exp(\A)=(1,d_2,\ldots,d_\ell)$ if $\A$ is not empty. 
\label{DH}
\end{prop}

Then we have the following algebro-geometric 
interpretation of $\chi(\A;t)$ by using $D_H(\A)$. 

\begin{theorem}[Theorem 4.1, \cite{MS}]
Assume that $\A$ is locally free, 
$H \in \A$ and 
let 
$\widetilde{D_H(\A)}$ denote the sheafification of $D_H(\A)$ onto 
$\P^{\ell-1} \simeq \mbox{Proj}(S)$. If $c_t(\widetilde{D_H(\A)})$ denotes the 
Chern polynomial of $\widetilde{D_H(\A)}$, then 
$$
t^{\ell-1} c_{1/t}(\widetilde{D_H(\A)})=\chi_0(\A;t).
$$
\label{MS}
\end{theorem}

\section{Proof of main results}

In this section we prove Theorem \ref{main}. First let us show the following 
easy results.

\begin{lemma}
Let  $H \in \A$ and $\A':=\A \setminus \{H\}$. Assume that $\A'$ is free, and 
$X \in L(\A^H)$. Then $\A_X':=\A_X \setminus \{H\}$ is free.
\label{LD0}
\end{lemma}

\noindent
\textbf{Proof}. 
First assume that 
$X \in L(\A')$ too. Then $(\A')_X=\A_X \setminus \{H\}$. Since $\A'$ is free, Theorem \ref{localization} shows that $\A_X'$ is free. 
Next assume that $X \not \in L(\A')$. Then there is $Y \in L(\A')$ such that 
$\codim Y+1=\codim X$ and that $X =Y 
\cap H$. Moreover, there are no $L \in \A'$ such that $L \cap Y=
X$. Then $\A_X \setminus \{H\}=(\A')_Y$, which is free by Theorem \ref{localization}.\owari 
\medskip

\begin{lemma}
If $\A$ and $\A':=\A \setminus \{H\}$ are both free, then $\A$ is divisional along $H$.
\label{LD}
\end{lemma}

\noindent
\textbf{Proof}.
Apply Theorems \ref{additiondeletion}, \ref{factorization}, and Lemma \ref{LD0}. 
\owari
\medskip

\begin{lemma}
Assume that $\A$ is divisional along $H \in \A$. Then for any $X \in L(\A^H)$, 
$\A_X$ is divisional along $H$. Moreover, if $\B_X$ denotes
the essential part of $\A_X$, then $\B_X$ is divisional along $H$ too.
\label{localdivisional}
\end{lemma}

\noindent
\textbf{Proof}. 
Immediate by Definitions \ref{divisional}, \ref{cha} and \ref{essentialpart}. \owari
\medskip

Now let us show the main result in this 
section, which is the key to prove Theorem \ref{main}.

\begin{theorem}
Assume that $\A$ is locally free, and $\A'$ is free. Then $\A$ is 
free if and only if $\chi(\A^H;t) \mid \chi(\A;t)$.
\label{maincor}
\end{theorem}

\noindent
\textbf{Proof}.
The statement is true if $\A'=\emptyset$. Assume not in the following. 
If $\A$ and $\A'$ are free, then Lemma \ref{LD} implies that 
$\A$ is divisional along $H$. So it suffices to show that 
$\A$ is free if $\A'$ is free, $\A$ is locally free and 
$\chi(\A^H;t) \mid \chi(\A;t)$. 
If $\exp(\A')=(1,d_2,\ldots,d_\ell)$ with 
$d_2 \le \cdots \le d_\ell$, then Theorem \ref{factorization} implies that 
$$
\chi_0(\A';t)=\prod_{i=2}^\ell (t-d_i).
$$
Since $\chi(\A^H;t) \mid 
\chi(\A;t)$, Theorem \ref{DR} shows that $\chi(\A^H;t) \mid 
\chi(\A';t)$ too, and thus 
$|\A'|-|\A^H|=d_j$ for some $j$. We may assume that $d_j<d_{j+1}$, or $j=\ell$. Let 
$I:=\{2,\ldots,\ell\} 
\setminus \{j\}$. 
Then by Theorems \ref{additiondeletion} and \ref{factorization}, 
\begin{eqnarray*}
\chi_0(\A;t)&=&(t-d_j-1)\prod_{i \in I} (t-d_i),\\
\chi_0(\A^H;t)&=&\prod_{i \in I} (t-d_i).
\end{eqnarray*}
Let $\theta_1=\theta_E,\theta_2,\ldots,\theta_\ell$ be a homogeneous 
basis for $D(\A')$ with $\deg \theta_i=d_i$. If $\theta_k 
\not \in D(\A)$ for some $k$ with $d_k=d_j$, then we can show that 
$\A$ is free. Namely, by Proposition \ref{B}, there is a homogeneous polynomial $B$ of 
degree $|\A'|-|\A^H|=d_j$ such that 
$$
\theta(\alpha_H) \in (\alpha_H,B)
$$
for all $\theta \in D(\A')$. So $\theta \in D(\A')$ is in $D(\A)$ if 
$\deg \theta <d_j$. Now assume that $\theta_k \not \in D(\A)$ for some $k$ with 
$d_k=d_j$. We may assume that $k=j$. Then $\theta_j
(\alpha_H)=B$ modulo $\alpha_H$. Let $\theta_i(\alpha_H)=a_i \alpha_H+b_i B$. Then replacing 
$\theta_i$ by $\theta_i-b_i \theta_j$, we may assume that 
$\theta_i \in D(\A)$ if $i \neq j$. Thus Theorem \ref{Saito} shows that 
$\theta_1,\theta_2,\ldots,\theta_{j-1},\alpha_H\theta_j,\theta_{j+1},\ldots,\theta_\ell$ form a 
basis for $D(\A)$. 

Assume not, i.e., 
$$
D':=\oplus_{i=2}^j S\theta_i \subset D(\A).
$$
By replacing $\theta_i \ (i=2,\ldots,j)$ by 
$\theta_i-(\theta_i(\alpha_H)/\alpha_H)\theta_E$, we may assume that 
$\theta_i \in D_H(\A)$ for $i=2,\ldots,j$. 
By composing the inclusion $D(\A) \subset D(\A')$ and 
the canonical projection $D(\A') \rightarrow D'$, we have a 
surjection 
$$
p:D(\A) \rightarrow D'
$$
which has a canonical section by the above. Hence 
$$
D(\A)=S\theta_E \oplus D' \oplus M
$$
for some $S$-subumodule $M$ of $D(\A)$. 
By Proposition \ref{DH}, 
$$
D_H(\A) \simeq D(\A)/S\theta_E=D' \oplus M.
$$
Since $\A$ is locally free, Theorem \ref{MS} implies that 
$$\chi_0(\A;t)=t^{\ell-1} c_{1/t}(\widetilde{D_H(\A)}).
$$
Since $D_H(\A) \simeq D' \oplus M$, it holds that 
$$
c_{1/t}(\widetilde{D_H(\A)})=
c_{1/t}(\widetilde{D'})c_{1/t}(\widetilde{M}).
$$
Noting that $D'=\oplus_{i=2}^j S[-d_i]$, 
$$
t^{\ell-1}c_{1/t}(\widetilde{D_H(\A)})=t^{\ell-j}c_{1/t}(\widetilde{M})\prod_{i=2}^j (t-d_i).
$$
Here 
$$
\rank(\widetilde{M})=
\rank(\widetilde{D_H(\A)})-
\rank(\widetilde{D'})=\ell-j.
$$
Since $\widetilde{D_H(\A)}$ is locally free, so is 
$\widetilde{M}$. Hence 
$\deg c_t(\widetilde{M})=\ell-j$. So 
$c(t):=c_{1/t}(\widetilde{M})t^{\ell-j} \in \Z[t]$. 
Therefore 
$$
\chi_0(\A;t)=(t-d_j-1)
\prod_{i \in I}^\ell (t-d_i)=
c(t)\prod_{i=2}^j(t-d_i).
$$
Assume that the multiplicity of $d_j$ in $\chi_0(\A';t)$ is $a \ge 1$. Then 
$d_j<d_{j+1}$ or $j=\ell$ implies that 
$C(t):=c(t)\prod_{i=2}^j(t-d_i)$ is divisible by 
$(t-d_j)^a$, but $\chi_0(\A;t)$ is not, a contradiction. \owari
\medskip

\noindent
\textbf{Proof of Theorem \ref{main}}. 
The statement is true if $\A'=\emptyset$. Assume not in the following. If $\A$ is free, then 
$\A$ is divisional along $H$ by Lemma \ref{LD}. So it suffices to show that $\A$ is free if $\A'$ is free, and 
$\A$ is divisional along $H$. 

We prove by induction on $\rank(\A) \ge 0$. We may assume that $\rank(\A)=\ell$ by 
Definition \ref{essentialpart}. The statement is true if $\ell \le 3$ by \cite{A}, or 
Theorems \ref{factorization} and \ref{division}. Assume that $\ell \ge 4$. By Theorem \ref{maincor}, it suffices to show that $\A$ is locally free. Since
$\A'$ is free, $\A_X=(\A')_X$ is free too if $X \not \subset H$. 
So let 
$X \in L_i(\A) \setminus \{0_V\}$ with $
X \subset H$. 
By Lemma \ref{LD0}, $\A_X \setminus \{H\}$ is free. 
Let $\B_X$ be the essential part of 
$\A_X$ as in Definition \ref{essentialpart}. In other words, $\A_X =\B_X \times \emptyset_{\ell-i}$, and $i =\rank(\B_X)< \rank (\A)=\ell$. Then by 
Proposition \ref{directfree}, $\B_X \setminus \{H\}$ is free too. 
Since $\B_X$ is of rank $i<\ell$, divisional along $H$ by Lemma \ref{localdivisional}, and 
$\B_X \setminus \{H\}$ is free, the induction hypothesis on $\ell$ shows that $\B_X$ is free. 
Thus $\A$ is locally free, and Theorem \ref{maincor} completes the proof. \owari
\medskip

\noindent
\textbf{Proof of Theorems \ref{additiondeletion2} and \ref{three}}. Immediate by 
Theorems \ref{additiondeletion}, \ref{main}, \ref{factorization}, \ref{division}, and \ref{deletioncomb}. \owari
\medskip

Finally we pose several conjectures related to Theorem \ref{main}. 

\begin{conj}
For $H \in \A$, let $\A':=\A \setminus \{H\}$. Then 
\begin{itemize}
\item[(1)] 
$\A$ is free if $\A'$ is free, and $\A$ is globally divisional along $H$. 
\item[(2)]
Assume that 
\begin{eqnarray*}
\chi_0(\A;t)&=&(t-d_2-1)\prod_{i=3}^\ell (t-d_i),\\
\chi_0(\A';t)&=&(t-d_2)\prod_{i=3}^\ell (t-d_i)
\end{eqnarray*}
for some integers $d_2,\ldots,d_\ell$. 
Then $\A$ and $\A'$ are both free.
\end{itemize}
\label{conj}
\end{conj}

Note that Conjecture \ref{conj} (1) and (2) are true 
if $\ell \le 3$ by Theorem \ref{main}, or Theorem \ref{division}.

\section{Additionally and stair-free arrangements of hyperplanes}

Theorem \ref{AFTC} is clear by Theorem \ref{main}. 
Let us show the relation among 
%
inductively, divisionally and additionally free arrangements. 

\begin{prop}
It holds that 
$$
\mathcal{AF} \supset \mathcal{IF} \subsetneq \mathcal{DF}.
$$
Moreover, $\mathcal{DF}  \not \subset \mathcal{AF}$.
\end{prop}

\noindent
\textbf{Proof}.
$\mathcal{IF} \subset \mathcal{AF}$ is clear by Definitions \ref{AF} and \ref{IF}. 
$\mathcal{IF} \subsetneq \mathcal{DF}$ follows by Theorem \ref{IFDF}. Moreover, 
by Lemma 3.13 in \cite{M}, 
we know that 
$\A(G_{31}) \in \mathcal{DF}  \setminus \mathcal{AF}$. Here 
$G_{31}$ is one of the finite unitary reflection groups 
classified and labeled by Shephard and Todd in \cite{ST}. So 
$\A(G_{31})$ is the corresponding unitary reflection arrangement in 
$\mathbb{C}^4$. 
\owari
\medskip

Now we have two class $\mathcal{DF}$ and $\mathcal{AF}$ of free arrangements which are larger than $\mathcal{IF}$. Since their definition is easy to deal with, we can join these two to 
obtain the new largest class of free arrangements in which 
Terao's conjecture is true.

\begin{define}[Stair-free arrangements]
The set $\mathcal{SF}_\ell$ consists of hyperplane arrangements $\A$ 
in $\K^\ell$ which satisfy the following; there are $H_1^i,\ldots,H_{j_i}^i, H_i\ 
(i=3,\ldots,\ell-1)$ in $L_{\ell-i}(\A)$ such that 
\begin{itemize}
\item[(1)]
$\A_\ell:=\A,\ 
\A_{i}:=(\A_{i+1} \setminus \{H_j^i\}_{j=1}^{j_i})^{H_i}$. So 
$H_j^i ,\ H_i\in \A_{i+1}$.
\item[(2)]
$\A_i^k:=\A_i \setminus \{ H_j^{i-1}\}_{j=1}^k$ is disivional 
along $H_{k+1}^{i-1}$ for all $i,k$.
\item[(3)]
$\chi(\A_i;t)$ divides $
\chi(\A_{i+1}^{j_i};t)$.
\end{itemize}
$\A \in \mathcal{SF}$ is called a \textbf{stair-free arrangement of hyperplanes}. 
\label{ADF}
\end{define}

By Theorems \ref{main} and \ref{division}, the following is clear.

\begin{theorem}
\begin{itemize}
\item[(1)]
$\A$ is free if $\A \in \mathcal{SF}$. 

\item[(2)]
It holds that 
$$\mathcal{IF}
 \subsetneq \mathcal{DF}\cup 
\mathcal{AF} \subset \mathcal{SF}.
$$

\item[(3)]
Terao's conjecture is true in the class of stair-free arrangements.
\end{itemize}
\label{ADFTC}
\end{theorem}

When you apply the division theorem, the arrangement moves vertically, i.e., the direction 
of dimensions, and the addition theorem horizontally, i.e., the direction of the cardinality of hyperplanes. If you connect these paths, then it looks like a stair. That is the reason of the name in 
Definition \ref{ADF}. A typical example of a stair-free arrangement is the Catalan and Shi  arrangements, which will be explained in \S5. In \cite{AT}, it is shown that the freeness of Catalan arrangement is combinatorial. Seeing the proof, \cite{AT} shows that Catalan arrangement is stair-free.

Related to these new classes, we have the following conjecture.

\begin{conj}
\begin{itemize}
\item[(1)]
$$
\mathcal{IF} \subsetneq \mathcal{AF}.
$$
\item[(2)]
$$
\mathcal{DF} \not \supset \mathcal{AF}.
$$
\item[(3)]
$$
\mathcal{DF} \cup \mathcal{AF} \subsetneq \mathcal{SF}.
$$
\end{itemize}
\end{conj}

\section{Applications to root systems}

In this section we apply Theorem \ref{main} to explicit examples, showing the combinatorial dependency of their freeness.

In this section let $\K=\R$, thus $V=\R^\ell$. Let $W$ an irreducible crystallographic Weyl group of rank $\ell$ acting on $V$. 
Let $\Phi$ be the corrsponding root system to $W$, $\Phi^+$ the set of positive roots, and 
$\Delta=\{\alpha_1,\ldots,\alpha_\ell\}$ the set of simple roots. 

A lot of arrangements related to $W$ are check to be free, and their freeness are often 
proved to be combinatorial. 
For example, if $\A_W$ denotes the set of all reflecting hyperplanes of reflections in $W$, then $\A_W \in \mathcal{IF}$, thus its freeness depends only on combinatorics. 
To review some of them, let us introduce some notation. 

For $\alpha \in \Phi^+$, let $H_\alpha$ denote the reflecting hyperplane of $\alpha$, and 
$$
H_\alpha^k:=\{\alpha=k\} \subset V
$$
be an affine hyperplane for $k \in \Z$. Let $cH_\alpha^k$ be the cone $\alpha=kz$ of $H_\alpha^k$ by the new coordinate in $z$, i.e., $cH_\alpha^k$ is a hyperplane in $\R^{\ell+1}=
\mbox{Spec}(S[z])$. 
Then for $m \in \Z_{>0}$, the \textbf{Catalan arrangement} $\Cat^m$ is defined by 
$$
\Cat^m:=\{z=0\} \cup \{cH_\alpha^k\}_{\alpha \in \Phi^+,-m \le k \le m},
$$
and the \textbf{Shi arrangement} $\Shi^m$ is defined by 
$$
\Shi^m:=\{z=0\} \cup \{cH_\alpha^k\}_{\alpha \in \Phi^+,-m+1 \le k \le m}.
$$
Both are free by \cite{Y1}, and their freeness depends only on the lattice by 
\cite{A2} for Shi arrangements, and \cite{AT} for Catalan arrangements. Also, we say that 
$I \subset \Phi^+$ is a \textbf{lower ideal} if 
$\alpha \in I,\ \beta \in \Phi^+,\ 
\alpha-\beta \in \sum_{i=1}^\ell \Z_{\ge 0} \alpha_i$ implies $\beta \in I$. Then it was shown that 
$\A_I:=\{H_\alpha \mid \alpha \in I\}$ is free in \cite{ABCHT}. Moreover, if we define the 
\textbf{ideal Shi arrangements} 
\begin{eqnarray*}
\Shi_{+I}^m:&=&\Shi^m \cup \{cH_\alpha^{-m}\mid 
\alpha \in I\},\\
\Shi_{-I}^m:&=&\Shi^m \setminus \{cH_\alpha^{m}\mid 
\alpha \in I\},
\end{eqnarray*}
then it is shown in Theorem 1.6 in \cite{AT0} that they are both free, and 
the freeness of $\Shi_{+m}^I$ depends only on combinatorics. So the left problems is that of $\Shi_{-I}^m$. Since ideal-Shi arrangements have a free filtration by 
Theorem 1.5 in \cite{AT0}, which is nothing but the additional filtration, 
Theorem \ref{AF} implies the following.

\begin{theorem}
It holds that
$$
\Shi_{\pm I}^m \in \mathcal{AF}.
$$
Thus, the freeness of $\Shi_{\pm I}^m$ depends only on 
$L(\Shi_{-I}^m)$ for all lower ideal $I \subset \Phi^+$ and $m \in \Z_{\ge 0}$. 
\label{idealShi}
\end{theorem}

Since $\Shi_{+\Phi^+}^m=\Cat^m$ and $\Shi_{\emptyset}^m=\Shi^m$, Theorem \ref{idealShi} gives a uniform proof of the combinatorial dependency of their freeness.

\end{document}